\documentclass[fleqn]{mat01}
\usepackage{times,mathtimy,amssymb,latexsym,amsbsy,rangecite,amsbsy}
\begin{document}

\setcounter{page}{201}
\firstpage{201}

\newtheorem{theore}{Theorem}
\renewcommand\thetheore{\arabic{section}.\arabic{theore}}
\newtheorem{theor}[theore]{\bf Theorem}
\newtheorem{rem}[theore]{Remark}
\newtheorem{propo}[theore]{\rm PROPOSITION}
\newtheorem{lem}{Lemma}
\newtheorem{definit}[theore]{\rm DEFINITION}
\newtheorem{assump}[theore]{Assumption}
\newtheorem{coro}[theore]{\rm COROLLARY}
\newtheorem{exampl}[theore]{Example}
\newtheorem{pot}[theore]{Proof of Theorem}
\newtheorem{case}{Case}

\def\corol{\trivlist \item[\hskip \labelsep{COROLLARY.}]}
\def\noteproof{\trivlist \item[\hskip \labelsep{\it Note added in Proof.}]}

\renewcommand{\theequation}{\thesection\arabic{equation}}

\title{A note on convex renorming and fragmentability}

\markboth{A K Mirmostafaee}{Convex renorming and fragmentability}

\author{A K MIRMOSTAFAEE}

\address{Department of Mathematics, Damghan University of
Sciences, P.O.~36715/364, Damghan, Iran\\
\noindent E-mail: kmirmostafaee@dubs.ac.ir}

\volume{115}

\mon{May}

\parts{2}

\pubyear{2005}

\Date{MS received 29 November 2004}

\begin{abstract}
Using the game approach to fragmentability, we give new and
simpler proofs of the following known results: (a)~If the Banach
space admits an equivalent Kadec norm, then its weak topology is
fragmented by a metric which is stronger than the norm topology.
(b)~If the Banach space admits an equivalent rotund norm, then its
weak topology is fragmented by a metric. (c)~If the Banach space
is weakly locally uniformly rotund, then its weak topology is
fragmented by a metric which is stronger than the norm topology.
\end{abstract}


\maketitle

\section{Introduction}

Let $(X, \tau)$ be a topological space and $\rho$ be a metric on
$X$. Given $\epsilon > 0$, a nonempty subset $A$ of $X$ is said to
be {\it fragmented by $\rho$ down to $\epsilon$} if each nonempty
subset of $A$ contains a nonempty $\tau$--relatively open subset
of $\rho$-diameter less than $\epsilon$. $A$ is called {\it
fragmented by $\rho$} if $A$ is fragmented by $\rho$ down to
$\epsilon$ for each $\epsilon > 0$. The set $A$ is said to be
$\sigma$-{\it fragmented by $\rho$} if for every $\epsilon > 0$,
$A$ can be expressed as $A = \cup_{n = 1}^{\infty} A_{n,
\epsilon}$ with each $A_{n, \epsilon}$ fragmented by $\rho$ down
to $\epsilon$.

The notion of fragmentability was originally introduced in
\cite{3} to investigate the existence of nice selections for upper
semicontinuous compact-valued mappings. The notion of
$\sigma$-fragmentability appeared in \cite{1} in order to study
Banach spaces, the weak topology of which is $\sigma$-fragmented
by the norm (such Banach spaces are said to be
$\sigma$-fragmentable). Since then, these two concepts have been
playing an important role in the study of the geometry of Banach
spaces.

Kenderov and Moors \cite{4} used the following topological game to
characterize fragmentability of a topological space $X$: Two
players $\Sigma$ and $\Omega$ alternatively select subsets of $X$.
$\Sigma$ starts the game by choosing some nonempty subset $A_{1}$
of $X$. Then $\Omega$ chooses some nonempty relatively open subset
$B_{1}$ of $A_{1}$. In general, if the selection $B_{n} \neq
\emptyset$ of the player $\Omega$ is already specified, the player
$\Sigma$ makes the next move by selecting an arbitrary nonempty
set $A_{n + 1}$ contained in $B_{n}$. Continuing the game the two
players generate a sequence of sets
\begin{equation*}
A_{1} \supset B_{1} \supset \cdots \supset A_{n} \supset B_{n}
\supset \cdots
\end{equation*}
which is called a play and is denoted by $p = (A_{i}, B_{i})_{i =
1}^{\infty}$. If
\begin{equation*}
p_{1} = (A_{1}), \ldots, p_{n} = (A_{1}, B_{1}, \ldots, A_{n})
\end{equation*}
are the first `$n$' move of some play (of the game), then $p_{n}$
is called the $n$th {\it partial play} of the game. The player
$\Omega$ is said to have won the play $p$ if $\cap_{i =
1}^{\infty} A_{i} = \cap_{i = 1}^{\infty} B_{i}$ contains at most
one point. Otherwise the player $\sum$ is said to be the winner in
this play. Under the term {\it strategy s for $\Omega$-player}, we
mean a rule by means of which the player $\Omega$ makes his/her
choices. More precisely, the strategy $s$ is a sequence of
mappings $s = \{s_{n}\}_{n \geq 1}$, which are defined inductively
as follows: $s_{1}$ assigns to each possible first move $A_{1}$ of
$\Sigma$-player a nonempty relatively open subset $B_{1} = s_{1}
(A_{1})$. Therefore, the domain of $s_{1}$ is the set of all
nonempty subsets of $X$ and $s_{1}$ assigns to each such an
element a nonempty relatively open subset of it. The domain of
$s_{2}$ consists of triples of the type $(A_{1}, B_{1}, A_{2})$,
where $A_{1}$ is from the domain of $s_{1}, B_{1} = s_{1}(A_{1})$
and $A_{2}$ is an arbitrary nonempty subset of $B_{1}$. $s_{2}$
assigns to such a triple a nonempty relatively open subset $B_{2}
= s_{2} (A_{1}, B_{1}, A_{2})$ of $A_{2}$. In general, the domain
of $s_{n + 1}$ consists of partial plays of the type
\begin{equation*}
(A_{1}, \ldots, A_{i}, B_{i}, A_{i + 1}, \ldots, A_{n + 1}),
\end{equation*}
where, for every $i \leq n, (A_{1}, \ldots , A_{i})$ is from the
domain of $s_{i}, B_{i} = s_{i} (A_{1}, \ldots, A_{i})$ and $A_{n
+ 1}$ is an arbitrary nonempty subset of $B_{n}$. To every element
from its domain $s_{n + 1}$ assigns a nonempty relatively open
subset $B_{n + 1}$ of $A_{n + 1}$.

A play $p = (A_{i}, B_{i})_{i \geq 1}$ is called an $s$-play if
$B_{i} = s_{i} (p_{i})$ for each $i \geq 1$. $s$ is called a {\it
winning strategy} for the player $\Omega$ if he/she wins every
$s$-play. If the space $X$ is fragmentable by a metric
$d(\cdot\,,\cdot)$, then $\Omega$ has an obvious winning strategy
$s$. Indeed, to each partial play $p_{n}$ this strategy puts into
correspondence some nonempty subset $B_{n} \subset A_{n}$ which is
relatively open in $A_{n}$ and has d-diameter less than $1/n$.
Clearly, the set $\cap_{i \geq 1} A_{i} = \cap_{i \geq 1} B_{i}$
has at most one point because it has zero d-diameter. Kenderov and
Moors have shown that the existence of a winning strategy for the
player $\Omega$ characterizes fragmentability, that is,

\begin{theor}[\cite{4}]
The topological space $X$ is fragmentable if and only if the
player $\Omega$ has a winning strategy.
\end{theor}

Of special interest is the case when the topology generated by the
fragmenting metric contains the original topology of the space (in
this case it is said that $X$ {\it is fragmented by a metric which
is stronger than its topology}).

\begin{theor}[\cite{4}]
The topological space $X$ is fragmentable by a metric stronger
than its topology if and only if the player $\Omega$ has a
strategy a such that{\rm ,} for every $s$-play $p = (A_{i},
B_{i})_{i \geq 1}$ the intersection $\cap_{i = 1}^{\infty} A_{i} =
\cap_{i = 1}^{\infty} B_{i}$ is either empty or contains just one
point $x_{0}$ and for every neighborhood $U$ of $x_{0}$ there
exists some $k$ such that $A_{i} \subset U$ for all $i > k$.
\end{theor}

This characterization of fragmentability has some applications
(see e.g. \cite{4,5,6}). In \cite{5}, it is shown that
fragmentability and $\sigma$-fragmentability of the weak topology
in a Banach space are related to each other in the following way:

\begin{theor}[(\cite{5}, Theorems~1.3, 1.4 and 2.1)] For a Banach
space $X$ the following are equivalent{\rm :}

\begin{enumerate}
\renewcommand\labelenumi{\rm (\roman{enumi})}
\leftskip .35pc
\item $(X,\,\hbox{weak})$ is $\sigma$-fragmented by the norm {\rm
(}i.e. $X$ is $\sigma$-fragmented{\rm );}

\item $(X,\,\hbox{weak})$ is fragmented by a metric which is
stronger than the weak topology{\rm ;}

\item $(X,\,\hbox{weak})$ is fragmented by a metric which is
stronger than the norm topology{\rm ;}

\item There exists a strategy $s$ for the player $\Omega$ in $(X,\,\hbox{weak})$
such that{\rm ,} for every $s$-play $p = (A_{i}, B_{i})_{i \geq
1}$ either $\cap_{i \geq 1} B_{i} = \emptyset$ or $\lim_{i
\rightarrow \infty}$ norm-diam $(B_{i}) = 0$.

\item There exists a strategy $s$ for the player $\Omega$ in $(X,\,\hbox{weak})$ such that{\rm ,} for every $s$-play $p = (A_{i},
B_{i})_{i \geq 1}$ either $\cap_{i \geq 1} B_{i} = \emptyset$ or
every sequence $\{x_{i}\}_{i \geq 1}$ with $x_{i} \in B_{i}, i
\geq 1$ has a weak cluster point.
\end{enumerate}
\end{theor}

Moreover, we have the following: The norm $\| \cdot\|$ of a Banach
space $X$ is said to be {\it Kadec} if the norm topology and the
weak topology coincide on the unit sphere $\{x \in X\hbox{:}\
\|x\| = 1\}$. In \cite{2}, it was shown that every Banach space
with Kadec norm is $\sigma$-fragmented. It follows that there
exists a strategy for the player $\Omega$ satisfying condition
(iv)~from the theorem of Kenderov and Moors. In the next section,
we will directly construct such a strategy (without using the
theorem of Kenderov and Moors).

The norm $\|\hbox{$\cdot$}\|$ of a Banach space $X$ is said to be {\it
rotund {\rm (}or strictly convex{\rm )}} if the unit sphere $\{x
\in X\,\hbox{:}\,\|x\| = 1\}$ does not contain nontrivial line
segments. Ribarska has shown in \cite{7} that the weak topology of
a rotund Banach space is fragmented by a metric. By the
abovementioned characterization of fragmentability it follows that
the player $\Omega$ has a winning strategy. In the next section we
will directly define such a strategy (without using the result of
Ribarska and the mentioned theorem of Kenderov and Moors).
Moreover, if the norm of $X$ is weakly locally uniformly rotund,
then the strategy we construct satisfies condition (v)~from the
above theorem of Kenderov and Moors. Recall that the Banach space
$X$ is called {\it locally uniformly rotund {\rm (}resp. weakly
locally uniformly rotund{\rm )}} if $\lim_{n \rightarrow \infty}
\|x_{n} - x\| = 0$ (resp. $\hbox{\it weak--}\lim(x_{n} - x) = 0$,
whenever $\lim_{n \rightarrow \infty} \|(x_{n} + x)/2\| = \lim_{n
\rightarrow \infty}\|x_{n}\| = \|x\|$.

\section{Description of the strategies}

\begin{lem}
Let $X$ be a Banach space with Kadec norm. Then{\rm ,} for every
$\epsilon > 0$ and $x \in X${\rm ,} there exists some positive
number $\alpha_{\epsilon, x}$ and a weakly open set $W_{\epsilon,
x} \ni x$ such that $\|y - x\| < \epsilon$ whenever $y \in
W_{\epsilon, x}$ and $|\|y\| - \|x\|| \leq \alpha_{\epsilon, x}$.
\end{lem}

\begin{proof}
If $x = 0$, it suffices to put $W_{\epsilon, x} = X$ and to take
as $\alpha_{\epsilon, x}$ any positive number smaller than
$\epsilon/2$. Suppose $x \neq 0$ and take a convex weakly open
neighborhood $G$ of $x$ such that the norm diameter of $G \cap
\{z\hbox{:}\ \|z\| = \|x\|\}$ is less than $\epsilon/2$. Define
$\alpha_{\epsilon, x} > 0$ to be smaller than $\epsilon/2, \|x\|$
and such that $\alpha_{\epsilon, x} B \subset (G - x)/2$ (as usual
$B$ stands for the closed unit ball of $X$). Put $W_{\epsilon, x}
:= x + (G - x)/2 = (x + G)/2$. Let $y \in W_{\epsilon,x}$ and
$|\|y\| - \|x\|| < \alpha_{\epsilon,x}$. Then we have
\begin{gather*}
(\|x\|/\|y\|) y = ((\|x\|/\|y\|)y - y) + y = (\|x\| -
\|y\|)y/\|y\| + y\\[.3pc]
\in |\|y\| - \|x\|| B + W_{\epsilon, x} \subset  \alpha_{\epsilon, x}
B + W_{\epsilon, x} \subset (G - x)/2 + (G + x)/2 = G.
\end{gather*}
Hence $\|(\|x\|/\|y\|)y - x)\| < \epsilon/2$. Finally we have
\begin{equation*}
\|y - x\| \leq \|y - (\|x\|/\|y\|)y\| + \|(\|x\|/\|y\|)y - x\| <
\alpha_{\epsilon, x} + \epsilon/2 < \epsilon.
\end{equation*}
$\left.\right.$\vspace{-2pc}

\hfill $\Box$
\end{proof}

We also need the following result:

\begin{lem}\hskip -.3pc {\rm (\cite{5}, Proposition~2.1).}\ \ 
If the closed unit ball $B$ of a Banach space $X$ admits a
strategy $s$ with the property {\rm (iv)} of Theorem~{\rm 1.3,}
then the whole space also admits such a strategy.
\end{lem}

\setcounter{theore}{0}
\begin{theor}[\!]
Let $X$ be a Banach space with Kadec norm. Then there exists a
strategy $s$ for the player $\Omega$ in {\rm (}$B${\rm ,} weak{\rm
)} such that{\rm ,} for every $s$-play $p = (A_{i}, B_{i})_{i \geq
1}$ either $\cap_{i \geq 1} B_{i} = \emptyset$ or $\lim_{i
\rightarrow \infty}$ norm-diam $(B_{i}) = 0$.
\end{theor}

\begin{proof}
Let $\|\hbox{$\cdot$}\|$ denote the Kadec norm on $X$ and $A_{1}$ be the
first choice of $\Sigma$-player. By Lemma~2, we may assume that
$A_{1} \subset B$, where $B$ denotes the closed unit ball of $X$.
Put
\begin{equation*}
\rho_{1} = \sup \{\|x\|\,\hbox{:}\,x \in A_{1}\}\quad
\hbox{and}\quad \epsilon_{1} = 1.
\end{equation*}
Two cases may happen.

\begin{enumerate}
\renewcommand\labelenumi{(\arabic{enumi})}
\leftskip .15pc

\item There is an element $x_{1} \in A_{1}$ such that
$\alpha_{\epsilon_{1}, x_{1}} + \|x_{1}\| > \rho_{1}$. Then we
take such a point $x_{1}$ and define $s_{1} (A_{1}) = B_{1} :=
W_{\epsilon_{1},x_{1}} \cap A_{1}\backslash(\|x_{1}\| -
\alpha_{\epsilon_{1}, x_{1}})B$ and $\epsilon_{2} :=
\epsilon_{1}/2$. Then for each $y \in B_{1}, \|y\| \leq \rho_{1} <
\alpha_{\epsilon_{1}, x_{1}} + \|x_{1}\|$ and $\|y\| \geq
\|x_{1}\| - \alpha_{\epsilon_{1}, x_{1}}$. Therefore, by Lemma~1,
$\|y - x_{1}\| < \epsilon_{1}$. Hence $\|\ \|-\hbox{diam}(B_{1}) <
2\epsilon_{1}$.

\item For every $x \in A_{1}, \alpha_{\epsilon_{1}, x} + \|x\|
\leq \rho_{1}$. Then,
\begin{equation*}
\hskip -1.2pc s_{1}(A_{1}) = B_{1} := A_{1}\backslash(1/2)
\rho_{1} B
\end{equation*}
and set $\epsilon_{2} = \epsilon_{1}$. Suppose the mappings
$(s_{i})_{i \leq n}$ participating in the definition of a strategy
for player $\Omega$ have already been defined. Let $(A_{i},
B_{i})_{1 \leq i \leq n}$ be a partial play which is generated by
the strategy mappings defined so far. This partial play is
accompanied by the numbers $\{\epsilon_{i}\}_{1 \leq i \leq n}$
and the points $x_{1}, \ldots, x_{n}$. If $A_{n + 1}$ is the next
move of the player $\Sigma$, we put
\begin{equation*}
\hskip -1.2pc \rho_{n + 1} = \sup \{\|x\|\,\hbox{:}\,x \in A_{n +
1}\}
\end{equation*}
and consider the following two possible cases:

\setcounter{enumi}{0}

\item There exists an element $x_{n + 1} \in A_{n + 1}$, such that
$\alpha_{\epsilon_{n + 1}, x_{n + 1}} + \|x_{n + 1}\| > \rho_{n +
1}$. In this case, we take such a point $x_{n + 1}$, define
\begin{align*}
\hskip -1.2pc s_{n + 1} (A_{1}, \ldots , A_{n + 1}) = B_{n + 1} :
= W_{\epsilon_{n + 1}, x_{n + 1}} \cap A_{n + 1} \backslash
(\|x_{n + 1}\| - \alpha_{\epsilon_{n}, x_{n + 1}})B
\end{align*}
and set $\epsilon_{n + 2} = \epsilon_{n + 1}/2$. As above one
shows that in this case $\|\ \|-\hbox{diam}(B_{n + 1}) < 2
\epsilon_{n + 1}$.

\item For every point $x \in A_{n + 1}, \alpha_{\epsilon_{n +
1},x} + \|x\| \leq \rho_{n + 1}$. In this case, we define
\begin{equation*}
\hskip -1.2pc s_{n + 1}(A_{1},\ldots,A_{n + 1}) = B_{n + 1} :=
A_{n + 1} \bigg\backslash \left(1 - \frac{1}{(n + 2)}\right)
\rho_{n + 1} B
\end{equation*}
and set $\epsilon_{n + 2} = \epsilon_{n + 1}$. In this way the
strategy $s = (s_{i})_{i \geq 1}$ for the $\Omega$-player is
already defined.
\end{enumerate}

Suppose $(A_{i}, B_{i})_{i \geq 1}$ is an $s$-play with $x \in
\cap_{n \geq 1} A_{n}$ and $\lim_{n \rightarrow \infty}
\|\hbox{$\cdot$}\|-\hbox{diam}(B_{n}) \neq 0$. Then there exists some
$\delta > 0$, such that $\|\hbox{$\cdot$}\|-\hbox{diam}(B_{n}) > \delta$
for each $n \in N$. This means that for all but finitely many $n$,
the case (2) happens and thus $\{\epsilon_{n}\}$ is eventually
constant: $\epsilon_{n} = \epsilon > 0$ for all $n > k$. Since $x
\in \cap_{n \geq 1} A_{n}$,
\begin{equation*}
\left(1 - \frac{1}{n}\right) \rho_{n} < \|x\| < \rho_{n},\quad
\hbox{for all}\quad n.
\end{equation*}
Let $\rho_{n}\!\hbox{$\searrow$}\rho$. Then the above inequality
shows that $\|x\| = \rho$. On the other hand, $\alpha_{\epsilon,
x} + \|x\| = \alpha_{\epsilon_{n}, x} + \|x\| \leq \rho_{n}$ for
$n > k$ which implies the contradiction $\alpha_{\epsilon, x} +
\|x\| = \|x\|$. \hfill $\Box$
\end{proof}

\begin{rem}
{\rm Lemma~1 directly implies that Banach spaces with Kadec norm
are $\sigma$-fragmentable. Actually, Theorem~2.3 of \cite{2}
indirectly implies that every Kadec renormable Banach space $X$
has a countable cover by sets of small local norm diameter, i.e.,
for each $\varepsilon > 0$, it is possible to write $X = \cup_{n
\in N} X_{n, \epsilon}$ such that for each $n \in N$ and $x \in
X_{n,\epsilon}$, there exists an open neighborhood $V_{x}$, of
$x$ such that the norm diameter of $V_{x} \cap X_{n, \epsilon}$
is less then $\epsilon$. Using Lemma~1, we can give another
proof of this result.}
\end{rem}

\begin{propo}$\left.\right.$\vspace{.5pc}

\noindent Let $X$ be a Banach space with Kadec norm. Then for
every $\epsilon > 0$ there exists a countable cover of $X, X =
\cup_{i \geq 0} X_{i}${\rm ,} such that{\rm ,} for every $x \in
X_{i}${\rm ,} there exists a weakly open neighborhood $W$ of $x$
such that $W \cap X_{i}$ is contained in $x + \epsilon B${\rm ,}
in particular the points of $X_{i}$ have weak neighborhoods with
norm-diameter smaller than $2\epsilon$.
\end{propo}

\begin{proof}
Given $\epsilon > 0$ consider, for $k = 1, 2, \ldots,$ and $n = 0,
1, 2, \ldots ,$ the sets $X_{kn} = \{ x \in X\,\hbox{:}\,
\alpha_{\epsilon, x} > 2/k$, and $n/k \leq \|x\| \leq (n +
1)/k\}$. Clearly, $X$ is covered by $X_{kn}$. Put $W :=
W_{\epsilon, x}$. By Lemma~1 the set $W \cap X_{kn}$ is contained
in $x + \epsilon B$.\hfill $\Box$
\end{proof}

\begin{theor}[\!]
Let $X$ be a Banach space.

\begin{enumerate}
\renewcommand\labelenumi{\rm (\alph{enumi})}
\leftskip .15pc
\item If the norm of $X$ is rotund{\rm ,} then {\rm (}$X${\rm ,}\,weak{\rm )} is fragmentable by a metric.

\item If the norm of $X$ is weakly locally uniformly rotund{\rm ,}
then {\rm (}$X${\rm ,}\,weak{\rm )} is fragmented by a metric which
is stronger than the norm topology.\vspace{-.5pc}
\end{enumerate}
\end{theor}

\begin{proof}
According to Theorems~1.2 and 1.3 and Lemma~2, it is enough to show
that in $(B,\,\hbox{\it weak})$ the player $\Omega$ has a winning
strategy $s$ such that, for every $s$-play $p = (A_{i}, B_{i})_{i
\geq 1}, \cap_{i \geq 1} B_{i}$ has at most one point and in case
(b) either $\cap_{i \geq 1} B_{i} = \emptyset$ or every sequence
$\{y_{n}\}, y_{n} \in B_{n}, n \geq 1$ is weakly convergent to the
element of $\cap_{i \geq 1} B_{i}$. Let $\|\ \|$ be the equivalent
norm on $X$ and $\Sigma$ start a game by choosing a nonempty
subset $A_{1}$ of $B$. Define
\begin{equation*}
\rho_{1} = \sup \{\|x\|\,:\,x \in A_{1}\}.
\end{equation*}
Choose an element $x_{1} \in A_{1}$ such that $\|x_{1}\| >
\rho_{1} - 1/2$ and find some $\mu_{1} \in X^{*}$ such that
$\|\mu_{1}\| = 1$ and $\mu_{1}(x_{1}) = \|x_{1}\|$. Define
\begin{equation*}
s_{1} (A_{1}) = B_{1} := \{x \in A_{1}\hbox{:}\ \mu_{1}(x) >
\rho_{1} - 1/2\}
\end{equation*}
as the first choice of $\Omega$-player. Then for each $x \in
B_{1}$, we have
\begin{equation*}
\rho_{1} - 1/2 < \mu_{1} (x) \leq \|x\| \leq \rho_{1}.
\end{equation*}

Suppose that the finite sequence $\{x_{k}\}_{k \leq n}$ of points
of $X, \{\mu_{k}\}_{k \leq n}$ of elements of $X^{*}$, and the
partial play $p_{n} = (A_{1}, \ldots, B_{n})$ have already been
specified so that for each $x \in B_{k}, k \leq n$ the inequality
\begin{equation*}
\rho_{k} - \frac{1}{k + 1} < \mu_{k} (x) < \|x\| \leq \rho_{k}
\end{equation*}
holds. Let $A_{n + 1}$ be the answer of $\Sigma$-player to
$p_{n}$. Put
\begin{equation*}
\rho_{n + 1} = \sup \{\|x\|\,\hbox{:}\,x \in A_{n + 1}\}
\end{equation*}
and find some $x_{n + 1} \in A_{n + 1}, \|x_{n + 1}\| > \rho_{n +
1} - \frac{1}{n + 2}$. Take some $\mu_{n + 1} \in X^{*}, \|\mu_{n
+ 1}\| = 1$ with $\mu_{n + 1} (x_{n + 1}) = \|x_{n + 1}\|$ and
define
\begin{align*}
s_{n + 1} (A_{1}, \ldots , A_{n + 1}) = B_{n + 1} = \left\lbrace x
\in A_{n + 1}\hbox{:}\ \mu_{n + 1}(x) > \rho_{n + 1} - \frac{1}{n
+ 2}\right\rbrace,
\end{align*}
as the next choice of the player $\Omega$. Clearly, for each $x
\in B_{n + 1}$, the inequality
\begin{equation*}
\rho_{n + 1} - \frac{1}{n + 2} < \mu_{n + 1} (x) < \|x\| \leq
\rho_{n + 1}
\end{equation*}
holds. Thus, by induction on $n$, we have shown that the
$\Omega$-player can choose sets of the form
\begin{equation*}
B_{n} = \left\{x \in A_{n}\,\hbox{:}\,\mu_{n} (x) > \rho_{n} -
\frac{1}{n + 1}\right\},
\end{equation*}
where $\|\mu_{n}\| = 1$ and $\rho_{n} = \sup \{\|x\|\,\hbox{:}\,x
\in A_{n}\}$ for each $n \in N$.

Let $\cap_{n \geq 1} B_{n} \neq \emptyset$ and $\mu$ be a
weak$^{*}$ cluster point of $\{\mu_{n}\}$. Then for each $x \in
\cap_{n \geq 1} B_{n}$, the inequality
\begin{equation*}
\rho_{n} - \frac{1}{n + 1} < \mu_{n} (x) < \|x\| \leq \rho_{n}
\end{equation*}
for each $n \in N$ implies that $\mu(x) = \|x\| = \rho$, where
$\rho$ is the limit of the decreasing sequence $\{\rho_{n}\}$. It
follows that for each $x, y \in \cap_{n \geq 1} B_{n}$, we have
$\mu(x) = \|x\| = \|y\| = \mu(y)$. Rotundity of $X$ implies that
$x = y$, thus, in this case, $\cap_{n \geq 1} B_{n}$ has at most
one point. In case (b), suppose that $x \in \cap_{n \geq 1}
B_{n}$. If $y_{n} \in B_{n}$, the inequality
\begin{equation*}
\rho_{n} - \frac{1}{n + 1} < \frac{1}{2} \mu_{n} (x + y_{n}) \leq
\frac{1}{2} \|x + y_{n}\| \leq \frac{1}{2} (\|x\| + \|y_{n}\|)
\leq \rho_{n}
\end{equation*}
shows that $\lim_{n \rightarrow \infty} \|(x + y_{n})/2\| =
\lim_{n \rightarrow \infty} \|y_{n}\| = \|x\| = \rho$. Since $(X,
\|\ \|)$ is weakly locally uniformly rotund, it follows that
$\lim_{n \rightarrow \infty} (x - y_{n}) = 0$. By Theorem~1.2, the
space is fragmented by a metric stronger than the weak topology.
This completes the proof.\hfill $\Box$
\end{proof}

\begin{rem}
{\rm It is well-known that locally uniformly rotund norms are
Kadec. Therefore statement (b) from the above theorem follows from
Theorem~2.1 as well. }
\end{rem}

\section*{Acknowledgement}

The author is indebted to P~Kenderov for the useful remarks while
this work was in progress.

\end{document}